\documentclass[12pt,reqno]{amsart}
\usepackage{amsmath,amsthm,amsfonts,amssymb,times}
\usepackage{verbatim}
\setbox0=\hbox{$+$}
\newdimen\plusheight
\plusheight=\ht0
\def\+{\;\lower\plusheight\hbox{$+$}\;}

\setbox0=\hbox{$-$}
\newdimen\minusheight
\minusheight=\ht0
\def\-{\;\lower\minusheight\hbox{$-$}\;}

\setbox0=\hbox{$\cdots$}
\newdimen\cdotsheight
\cdotsheight=\plusheight
\def\cds{\lower\cdotsheight\hbox{$\cdots$}}

\renewcommand{\Im}{\operatorname{Im}}

\makeatletter
\def\leqalignno#1{\displ@y \tabskip\z@ plus\@ne fil
  \halign to\displaywidth{\hfil$\@lign\displaystyle{##}$\tabskip\z@skip
    $\@lign\displaystyle{{}##}$\hfil\tabskip\z@ plus\@ne fil
    &\kern-\displaywidth\rlap{$\@lign\hbox{\rm##}$}\tabskip\displaywidth\crcr
    #1\crcr}}
\makeatother

\newcommand{\eb}{\begin{equation}}
\newcommand{\ee}{\end{equation}}
\renewcommand{\Im}{\operatorname{Im}}
\newcommand{\df}{\dfrac}
\newcommand{\tf}{\tfrac}

\renewcommand{\Im}{\operatorname{Im}}

\renewcommand{\Im}{\text{Im}}

\renewcommand{\(}{\left\(}
\renewcommand{\)}{\right\)}
\renewcommand{\[}{\left\[}
\renewcommand{\]}{\right\]}
\renewcommand{\i}{\infty}

\numberwithin{equation}{section}
 \theoremstyle{plain}
\newtheorem{theorem}{Theorem}[section]

\newtheorem{entry}[theorem]{Entry}

\allowdisplaybreaks

\begin{document}

\title{Explicit Values for Ramanujan's Theta Function $\varphi(q)$}

\author{Bruce C.~Berndt and \"{O}rs Reb\'{a}k}
\address{Department of Mathematics, University of Illinois, 1409 West Green
Street, Urbana, IL 61801, USA}\email{berndt@illinois.edu}
\address{Department of Mathematics and Statistics, University of Troms\o{} -- The Arctic University of Norway, 9037 Troms\o{}, Norway}
\email{ors.rebak@uit.no}

\begin{abstract}
This paper provides a survey of particular values of Ramanujan's theta function $\varphi(q)=\sum_{n=-\infty}^{\infty}q^{n^2}$, when $q=e^{-\pi\sqrt{n}}$, 
where $n$ is a positive rational number. First, descriptions of the tools used to evaluate theta functions are given.  Second, classical values are briefly discussed. 
Third, certain values due to Ramanujan and later authors are given.  Fourth, the methods that are used to determine these values are described.  
Lastly, an incomplete evaluation found in Ramanujan's lost notebook, but now completed and proved, is discussed with a sketch of its proof.
\end{abstract}


\maketitle

\centerline{\emph{Dedicated to the memory of Srinivasa Ramanujan}}

\section{Introduction}  Ramanujan loved to find closed-form evaluations of many items, e.g., definite integrals, infinite series, infinite products, class invariants, singular moduli, and theta functions.  Theta functions were at the epicenter of a significant portion of his research. In his final work on mock theta functions, the behavior of theta functions near their boundary of convergence on the unit circle was perhaps his chief motivating factor.  Not only did Ramanujan enjoy calculating special values of individual theta functions, but he also had a marvellous  insight for finding certain quotients of theta functions that yield elegant evaluations.    The purpose of this paper is to provide a survey of explicit values of perhaps the most important theta function, $\varphi(q)$, in Ramanujan's notation, which we define below.  In the final section, special attention is given to a mysterious, incomplete identity found in Ramanujan's lost notebook \cite[p.~206]{lnb}, which George Andrews and the first author left unfinished in their second volume on the lost notebook \cite[p.~181]{abII}.  Ramanujan's enigmatic entry has now been completed and proved by the second author \cite{rebak}.

\section{Ramanujan's theta functions}
Ramanujan's most general theta function $f(a,b)$ is defined by \cite[Volume~2, p.~197]{nb},\linebreak \cite[p.~34]{III}
\begin{equation}\label{ab}
f(a,b):=\sum_{n=-\infty}^{\infty}a^{n(n+1)/2}b^{n(n-1)/2}, \qquad |ab| <1.
\end{equation}
In classical notation,  $a=qe^{2iz}$ and $b=qe^{-2iz}$, where $|q|<1$, $z \in \mathbb{C}$, and $\Im(z) > 0$.
This definition of a theta function apparently originates with Ramanujan, i.e., to the best of our knowledge, no previous researcher had defined a general theta function by \eqref{ab}. For his purposes, the notation \eqref{ab} was far more advantageous and easier to use than the classical notation.  In particular, note the symmetry in $a$ and $b$ in \eqref{ab}, i.e., $f(a,b)=f(b,a)$.  Ramanujan loved symmetry.  Often, in expressing a function, identity, or theorem, if it were possible to state it symmetrically, he would do so.
 The symmetry reflected in the definition of $f(a,b)$ is inherited by its representation by the Jacobi triple product identity,
perhaps the most useful property of theta functions, given by \cite[Volume~2, p.~197]{nb}, \cite[p.~35, Entry~19]{III}
\begin{equation}\label{fabjtp}
f(a,b) = (-a;ab)_{\i}(-b;ab)_{\i}(ab;ab)_{\i}, \qquad |ab|<1,
\end{equation}
where
\begin{equation*}
(a;q)_{\i}:=\lim_{n\to\i}(a;q)_n,  \qquad |q|<1,
\end{equation*}
and
\begin{equation*}
(a;q)_n :=\prod_{k=0}^{n-1}(1-aq^k), \qquad n\geq1, \qquad (a;q)_0:=1.
\end{equation*}

In Ramanujan's notation, the three most important special cases of $f(a,b)$ in their series and product representations from \eqref{ab} and \eqref{fabjtp}, respectively, are defined by
\begin{align}
 \varphi(q) :=& f(q,q) =\sum_{n=-\i}^{\i}q^{n^2} =
 (-q;q^2)_{\i}^2(q^2;q^2)_{\i},\label{phi}\\
\psi(q) :=&f(q,q^3) =
\sum_{n=0}^{\i}q^{n(n+1)/2}=\df{(q^2;q^2)_{\i}}{(q;q^2)_{\i}},\label{psi}\\
\intertext{and}
 f(-q) :=& f(-q,-q^2) =\sum_{n=-\i}^{\i}(-1)^n
q^{n(3n-1)/2}=(q;q)_{\i}.\label{f}
\end{align}

In this paper, we focus on explicit values of $\varphi(q)$ of the form
\begin{equation}\label{n}
\varphi(e^{-\pi\sqrt{n}}),
\end{equation}
where $n$ is a positive rational number. In particular, concentration is given in those cases when $n$ is the square of a positive  integer. The representation \eqref{n} appears naturally in the theory of class invariants and singular moduli, which, along with modular equations, provide the most useful known tools for determining exact values of theta functions.  Ramanujan apparently used these connections to calculate several original values of  $\varphi(e^{-n\pi})$.

\section{Background Needed for the Determination of Values for $\varphi(q)$}
Recall that the ordinary hypergeometric function ${_2F_1}$ is defined for $|z|<1$ by
\begin{equation*}
{_2F_1}(a,b;c;z):=\sum_{n=0}^{\i}\df{(a)_n(b)_n}{(c)_nn!}z^n,
\end{equation*}
 where $(a)_0:=1$ and
$(a)_n:=a(a+1)(a+2)\cdots(a+n-1), n\geq1.$
 We use Ramanujan's notation to state one of the fundamental results in the classical theory of elliptic and theta functions, namely
 \cite[Volume~2, p.~207]{nb}, \cite[p.~101, Entry~6]{III},
    \begin{equation}\label{zF}
 z:= {_2F_1}\left(\tf12,\tf12;1;x\right)=\varphi^2(q),
 \end{equation}
 when
 \begin{equation}\label{q}
  q:=\exp\left(-\pi\,\df{{_2F_1}(\tf12,\tf12;1;1-x)}{{_2F_1}(\tf12,\tf12;1;x)}\right)=:e^{-y}.
  \end{equation}

With $x, y,$ and $z$ related as in \eqref{zF} and \eqref{q}, consider an equation in the form
\begin{equation}\label{zero}
\Omega(x,e^{-y},z)=0.
\end{equation}
This then implies an equation of the form
\begin{equation}\label{dup}
\Omega\left(\left\{\df{1-\sqrt{1-x}}{1+\sqrt{1-x}}\right\}^2, e^{-2y},\df12z\left(1+\sqrt{1-x}\right)\right)=0,
\end{equation}
which we call \emph{obtaining a formula by duplication}. The equation \eqref{zero} also implies that
\begin{equation}\label{dim}
\Omega\left(\df{4\sqrt{x}}{(1+\sqrt{x})^2},e^{-y/2},z(1+\sqrt{x})\right)=0,
\end{equation}
which we call \emph{obtaining a formula by dimidiation}.  Lastly, \eqref{zero} implies that
\begin{equation}\label{sign}
\Omega\left(\df{x}{x-1},-e^{-y},z\sqrt{1-x}\right)=0,
\end{equation}
which we call \emph{obtaining a formula by change of sign}.
Proofs for all three processes can be found in \cite[pp.~125, 126]{III}.

  Modular equations, singular moduli, and class invariants are the keys to determining specific values of $\varphi(e^{-\pi\sqrt{n}})$.  To define a modular equation, we first need to define the complete elliptic integral of the first kind $K(k)$, namely,
 \begin{equation}\label{ellipticintegral}
 K(k):=\int_0^{\pi/2}\df{dt}{\sqrt{1-k^2\sin^2t}}=\dfrac{\pi}{2}\,{_2F_1}(\tf12,\tf12;1;k^2), \qquad |k|<1.
 \end{equation}
 The number $k$ is called the \emph{modulus}, and $k^\prime:=\sqrt{1-k^2}$ is the \emph{complementary modulus}.

  Let $K, K^\prime, L$, and $L^\prime$ be complete elliptic integrals of the first kind associated with the moduli $k, k^{\prime}, \ell$, and $\ell^{\prime}$, respectively.  In his notebooks \cite{nb}, Ramanujan always used the notations $\alpha =k^2$ and $\beta=\ell^2$, and so we shall also  do this in the remainder of this article.  Suppose that for some positive integer $n$, the equality
 \begin{equation}\label{1}
 n\df{{_2F_1}(\tf12,\tf12;1;1-\alpha)}{{_2F_1}(\tf12,\tf12;1;\alpha)}
 =\df{{_2F_1}(\tf12,\tf12;1;1-\beta)}{{_2F_1}(\tf12,\tf12;1;\beta)}
 \end{equation}
 holds. Then a \emph{modular equation of degree $n$} is a relation between $\alpha$ and $\beta$ that is induced by \eqref{1}.
 The \emph{multiplier} $m$ for a modular equation of degree $n$ is defined by \cite[p.~230]{III}
 \begin{equation}\label{m}
  m:=\df{\varphi^2(q)}{\varphi^2(q^n)}=\df{{_2F_1}(\tf12,\tf12;1;\alpha)}{{_2F_1}(\tf12,\tf12;1;\beta)},
  \end{equation}
  by \eqref{zF}, where $x$ has now been replaced by $\alpha$.  From a modular equation of degree $n$ involving $\alpha, \beta,$ and $m$, we can obtain another modular equation, called the \emph{reciprocal modular equation}, which we now define \cite[Volume~2, p.~228]{nb}, \cite[p.~216, Entry~24(v)]{III}.

  \begin{entry}\label{reciprocal}
  If we replace $\alpha$ by $1-\beta$, $\beta$ by $1-\alpha$, and $m$ by $n/m$, where $n$ is the degree of the modular equation, we obtain a modular equation of the same degree.
  \end{entry}

  Let $n$ be a positive rational number.  Referring to \eqref{zF} and \eqref{ellipticintegral}, we define $\alpha_n$ by
  \begin{equation}\label{singular}
  \varphi^2(e^{-\pi\sqrt{n}})={_2F_1}(\tf12,\tf12;1;\alpha_n)
  =\df{2}{\pi}K(\sqrt{\alpha_n}).
  \end{equation}
 Then $\sqrt{\alpha_n}$ is called a \emph{singular modulus}.  The equation \eqref{singular} shows that the values of theta functions, hypergeometric functions, and complete elliptic integrals are intimately related, i.e.,  the evaluation of any one of these three quantities in \eqref{singular} yields a value for each of the other two objects. In the literature, perhaps more attention has been given to the evaluation of $K(\sqrt{\alpha_n})$. The techniques that are used typically express the values of $K(\sqrt{\alpha_n})$ in terms of gamma functions.  Moreover, Selberg and Chowla \cite{chowlaselberg} showed that for any singular modulus $\sqrt{\alpha_n}$, $K(\sqrt{\alpha_n})$ can be expressed in terms of gamma functions.  For specific evaluations, see the papers by J.~M.~Borwein and I.~J.~Zucker \cite{borweinzucker} and Zucker \cite{zucker}. A complete list of the values of $\varphi(e^{-\pi\sqrt{n}})$, $1\leq n\leq 16$, can be found in the  well-known treatise \cite[p.~298]{borwein} by J.~M.~Borwein and P.~B.~Borwein.

 Lastly, we define a \emph{class invariant}.  Let
 \begin{equation}\label{chi}
 \chi(q):=(-q;q^2)_{\infty}=2^{1/6}\left\{\alpha(1-\alpha)/q\right\}^{-1/24}, \qquad |q|<1,
 \end{equation}
 where in the latter representation, $q$ is given by \eqref{q} (with $x$ replaced by $\alpha$), and where a proof of the latter representation of $\chi(q)$ can be found in \cite[p.~124]{III}.
 If $q=e^{-\pi\sqrt{n}}$, where $n$ is a positive rational number, then the \emph{class invariant} $G_n$ is defined by
 \begin{equation}\label{G}
 G_n:=2^{-1/4}q^{-1/24}\chi(q).
 \end{equation}
 Using the latter representation for $\chi$ in \eqref{chi} and \eqref{G}, we deduce that
 \begin{equation}\label{Gn}
 G_n=\left\{4\alpha(1-\alpha)\right\}^{-1/24}.
 \end{equation}
 If $q=e^{-\pi}$ and  $\beta$ has degree $n$ over $\alpha$, it follows from \eqref{Gn} that
 \begin{equation}\label{Gn2}
 G_{n^2}=\left\{4\beta(1-\beta)\right\}^{-1/24}.
 \end{equation}

 To explicitly determine a value of $\varphi(e^{-n\pi })$ for a certain positive integer $n$, we choose an appropriate modular equation[s] of degree $n$ that frequently contains the multiplier $m$, given by \eqref{m}.  We now realize that we are free to choose any convenient value of $\alpha$, $0<\alpha<1$.  We thus will obtain an equation[s] involving $m$ and $\beta$.  Our goal is to express our equation[s] in terms of $m$ and a class invariant $G_{n^2}$, whose value is known.  Amazingly, Ramanujan calculated a total of 116 different class invariants.  See a complete table of Ramanujan's class invariants in \cite[pp.~189--204]{V}.

\section{Classical Values}
 The following values for $\varphi(e^{-\pi})$, $\varphi(e^{-\pi\sqrt2})$, and $\varphi(e^{-2\pi})$, are classical and were also discovered by Ramanujan \cite[Volume~2, p.~207]{nb}, \cite[pp.~103, 104, Entry~6]{III}:
\begin{equation}\label{classical}
\varphi(e^{-\pi})=\df{\pi^{1/4}}{\Gamma(\tf34)}, \quad \varphi(e^{-\pi\sqrt2})=\df{\Gamma(\tf98)}{\Gamma(\tf54)}\sqrt{\df{\Gamma(\tf14)}{2^{1/4}\pi}},\quad \text{and} \quad \varphi(e^{-2\pi})=\df{\sqrt{2+\sqrt2}}{2}\df{\pi^{1/4}}{\Gamma(\tf34)}.
\end{equation}
In fact, much more is true.  By the processes of \emph{duplication}, \emph{dimidiation}, and \emph{change of sign},  Ramanujan expressed  $$\varphi(-e^{-y}), \varphi(e^{-2y}), \varphi(-e^{-2y}), \varphi(e^{-4y}), \varphi(e^{-y/2}), \varphi(-e^{-y/2}), \varphi(e^{-y/4}),\quad \text{and } \quad \varphi(-e^{-y/4})$$ in  terms of the modulus $\sqrt{\alpha}$ and $\sqrt{z}$, which is a factor in each of these expressions
 \cite[Volume~2, p.~210]{nb}, \cite[p.~122, Entry~10]{III}.  For example,
\begin{equation}\label{zalpha}
\varphi(e^{-4y})=\tf12\sqrt{z}\left(1+(1-\alpha)^{1/4}\right).
\end{equation}
 Thus, beginning with  $\varphi(e^{-\pi})$ and  repeating the aforementioned three processes, we can obtain an infinite family of evaluations that can be added to those in \eqref{classical}.

  For the remainder of this paper, concentration is given to
$\varphi(e^{-n\pi}), n\geq3$.

\section{Values of $\varphi(e^{-n\pi })$ Found by Ramanujan;\linebreak Work of Heng Huat Chan and the First Author}

While in England, Ramanujan submitted a problem to the \emph{Journal of the Indian Mathematical Society} \cite{629}, wherein the second part   is equivalent to establishing the value
\begin{equation}\label{5}
\varphi(e^{-5\pi})=\df{\varphi(e^{-\pi})}{\sqrt{5\sqrt5-10}},
\end{equation}
which can be found in both Ramanujan's first \cite[Volume~1, pp.~285]{nb}, \cite[p.~327]{V} and second notebooks \cite[Volume~2, p.~227]{nb}, \cite[pp.~209, 210, Entry~23]{III}.  There were three claimants for a solution to Ramanujan's problem, one of which was incorrect.  For the first part of his question, Ramanujan asked readers to prove that, for $|x|<1$,
$$\dfrac12+\sum_{n=1}^{\infty}e^{-\pi n^2x}\cos\left\{\pi n^2\sqrt{1-x^2}\right\}=
\df{\sqrt2+\sqrt{1+x}}{\sqrt{1-x}}\sum_{n=1}^{\infty}e^{-\pi n^2x}\sin\left\{\pi n^2\sqrt{1-x^2}\right\}.$$
For further discussion,  see \cite[pp.~32, 33]{lange}.

In addition to the value of $\varphi(e^{-5\pi})$, values of $\varphi(e^{-3\pi})$, $\varphi(e^{-7\pi})$, $\varphi(e^{-9\pi})$, and $\varphi(e^{-45\pi})$ were also recorded by Ramanujan in his first notebook \cite[Volume~1, pp.~284, 297, 287, 312]{nb}, \cite[pp.~327, 328]{V}.
They were first proved in print by Heng Huat Chan and the first author \cite{cb}.  We record these four values:
\begin{align}
\df{\varphi(e^{-3\pi})}{\varphi(e^{-\pi})}&=\df{1}{\sqrt[4]{6\sqrt{3}-9}},\label{3}\\
 \df{\varphi^2(e^{-7\pi})}{\varphi^2(e^{-\pi})}&=
\df{\sqrt{13+\sqrt7}+\sqrt{7+3\sqrt7}}{14}(28)^{1/8},\label{77}\\
\label{999}
\df{\varphi(e^{-9\pi})}{\varphi(e^{-\pi})}&=\df{1+\sqrt[3]{2(\sqrt3+1)}}{3},\\
\intertext{and}
\df{\varphi(e^{-45\pi})}{\varphi(e^{-\pi})}&
=\df{3+\sqrt5+\left(\sqrt3+\sqrt5+(60)^{1/4}\right)\sqrt[3]{2+\sqrt3}}{3\sqrt{10+10\sqrt5}}.
\notag
\end{align}

 We provide a proof of only the evaluation \eqref{3}; it is taken from \cite{cb} and \cite[pp.~329, 330]{V}.

\begin{proof}
 If $\beta$ has degree $3$ over $\alpha$, then one of Ramanujan's 15 modular equations of degree $3$ and its reciprocal modular equation are  given by \cite[Volume~2, p.~230]{nb}, \cite[p.~230,\linebreak Entry~5(vii)]{III}
\begin{align}
m^2&=\left(\df{\beta}{\alpha}\right)^{1/2} +\left(\df{1-\beta}{1-\alpha}\right)^{1/2} -\left(\df{\beta(1-\beta)}{\alpha(1-\alpha)}\right)^{1/2}\label{first}\\
\intertext{and}
\df{9}{m^2}&=\left(\df{\alpha}{\beta}\right)^{1/2} +\left(\df{1-\alpha}{1-\beta}\right)^{1/2} -\left(\df{\alpha(1-\alpha)}{\beta(1-\beta)}\right)^{1/2}.\label{second}
\end{align}
Set $\alpha=1/2$ in \eqref{first} and \eqref{second}.  Next, multiply both sides of \eqref{second} by $2\{\beta(1-\beta)\}^{1/2}$ and then subtract the result from \eqref{first}.  This gives
\begin{equation*}
m^2-2\{\beta(1-\beta)\}^{1/2}\df{9}{m^2}=1-2\{\beta(1-\beta)\}^{1/2},
\end{equation*}
which, with the use of \eqref{Gn2} with $n=3$, yields
\begin{equation}\label{GG}
m^2-\df{9}{m^2G_9^{12}}=1-\df{1}{G_9^{12}}.
\end{equation}
Multiply both sides of \eqref{GG} by $G_9^6$ and use the value
\cite[p.~189]{V}
\begin{equation}\label{G9}
G_9=\left(\df{1+\sqrt3}{\sqrt2}\right)^{1/3}
\end{equation}
to arrive at
\begin{equation*}
(G_9^3m)^2-\df{9}{(G_9^3m)^2}=G_9^6-G_9^{-6}=2\sqrt3.
\end{equation*}
Hence, $(G_9^3m)^2=3\sqrt3$, and, by \eqref{G9}, $m^2=6\sqrt3-9$.
If we now appeal to \eqref{m}, we complete the proof of  \eqref{3}.
\end{proof}

In \cite{cb}, the first author and Chan also established explicit values for $\varphi(e^{-13\pi})$, $\varphi(e^{-27\pi})$, and $\varphi(e^{-63\pi})$. Values of the associated hypergeometric series were also derived.  As Ramanujan undoubtedly did to determine his values, in their proofs, these authors also used Ramanujan's modular equations and class invariants. Next, we provide the three above-mentioned values for $\varphi(q)$. First, define

  \begin{align*}
  G:=G_{169} =& \df{1}{3}\left(\sqrt{13}+2+\left(\df{13+3\sqrt{13}}{2}\right)^{1/3}\right.\\
 &\left.\times\left\{\left(\df{11+\sqrt{13}}{2}+3\sqrt3\right)^{1/3}
 +\left(\df{11+\sqrt{13}}{2}-3\sqrt3\right)^{1/3}\right\}\right)
  \end{align*}
  and
  \begin{equation*}
  a:=(G-G^{-1})^3+7(G-G^{-1}).
  \end{equation*}
  Then
  \begin{equation*}
  \df{\varphi(e^{-13\pi})}{\varphi(e^{-\pi})}=\left(G^{-3}\left(\df{a+\sqrt{a^2+52}}{2}\right)\right)^{-1/2}.
  \end{equation*}

  Next,
  \begin{equation}\label{10}
  \df{\varphi(e^{-27\pi})}{\varphi(e^{-3\pi})}
  =\df13\left(1+(\sqrt3-1)\left(\df{\sqrt[3]{2(\sqrt3+1)}+1}{\sqrt[3]{2(\sqrt3-1)}-1}\right)^{1/3}\right).
  \end{equation}
  By combining the identity \eqref{10} with the value \eqref{3}, we obtain the value of $\varphi(e^{-27\pi})$.

  Lastly,
  \begin{align}\label{11}
  \df{\varphi(e^{-63\pi})}{\varphi(e^{-7\pi})}=
  \df13&\left(1+\left(\df{\sqrt{4+\sqrt7}-7^{1/4}}{2}\right)^{3}\sqrt{\sqrt3+\sqrt7}(2+\sqrt3)^{1/6}\right.\notag\\
  &\left.\times\sqrt{\df{2+\sqrt7+\sqrt{7+4\sqrt7}}{2}}\sqrt{\df{\sqrt{3+\sqrt7}+(6\sqrt7)^{1/4}}
  {\sqrt{3+\sqrt7}-(6\sqrt7)^{1/4}}}\,\right).
    \end{align}
  Combining the evaluations \eqref{11} and \eqref{77}, we determine the value of $\varphi(e^{-63\pi})$.

\section{Values of $\varphi(e^{-n\pi})$, Modular Equations; Contributions of Jinhee Yi}
 In her paper \cite{yi}, Jinhee Yi established several new values for $\varphi(e^{-n\pi})$. In subsequent papers \cite{yi2} and \cite{yi3} with her colleagues, further new values were derived.  We next briefly describe her work and offer a few of her new values of $\varphi(e^{-n\pi})$.

 For any positive real numbers $n$ and $k$, define
 \begin{equation}\label{h}
 h_{k,n}:=\df{\varphi(e^{-\pi\sqrt{n/k}})}{k^{1/4}\varphi(e^{-\pi\sqrt{nk}})} \qquad \text{and}
 \qquad h_{k,n}^{\prime}:=\df{\varphi(-e^{-2\pi\sqrt{n/k}})}{k^{1/4}\varphi(-e^{-2\pi\sqrt{nk}})}.
 \end{equation}
 Two similar quotients involving the Dedekind eta-function (or $f(-q)$ defined in \eqref{f}) may be defined.  In an elementary way, Yi derived several relations among these four quotients.  The following theorem is an example \cite[p.~387]{yi}.

 \begin{theorem} For all positive real numbers $k, a, b, c$, and $d$, with $ab=cd$,
 $$ h_{a,b}h_{kc,kd}= h_{ka,kb}h_{c,d}.$$
 \end{theorem}

  She next stated or derived four modular equations (of degrees 4, 9, 15, 15).  We give one of the modular equations of degree 15 \cite[p.~391]{yi}, \cite[p.~235]{IV}.

  \begin{theorem} Let
  $$ P=\df{\varphi(q)}{\varphi(q^5)} \qquad \text{and}\qquad Q=\df{\varphi(q^3)}{\varphi(q^{15})}. $$
  Then
  $$ PQ+\df{5}{PQ}=\left(\df{Q}{P}\right)^2+3 \left(\df{Q}{P}\right)+3 \left(\df{P}{Q}\right) -\left(\df{P}{Q}\right)^2.$$
  \end{theorem}

  Employing quotients of theta functions, in particular, those in \eqref{h}, in the aforementioned four modular equations, Yi found several new values of $\varphi(e^{-\pi\sqrt{n}})$ \cite[pp.~391, 394, 396, 398--400]{yi}.  We offer only a small sampling:
  \begin{gather*}
  \df{\varphi(e^{-\sqrt3\pi})}{3^{1/4}\varphi(e^{-3\sqrt3\pi})}=\df{1}{\sqrt3}\left(1-\sqrt[3]{2}+\sqrt[3]{4}\right), \qquad
   \df{\varphi(e^{-\sqrt{5/3}\,\pi})}{3^{1/4}\varphi(e^{-\sqrt{15}\,\pi})}=\df{\sqrt{\sqrt5-1}}{\sqrt2},\\
   \intertext{and}
   \df{\varphi(-e^{-6\pi})}{\varphi(e^{-\pi})}
   =\df{(1+\sqrt3+\sqrt2\,\sqrt[4]{3^3})^{1/3}}{2^{11/24}3^{3/8}(\sqrt3-1)^{1/6}}.
  \end{gather*}

  This study continues in \cite{yi2} and \cite{yi3}, where further modular equations of `small' degree are used in conjunction with quotients of theta functions, including \eqref{h}.  An example from  \cite[p.~1325]{yi2} follows:
  \begin{equation*}
  \df{\varphi(e^{-2\pi/\sqrt5})}{5^{1/4}\varphi(e^{-2\sqrt5\pi})}=\df{2\sqrt{2a}}{(3+\sqrt2+\sqrt5+\sqrt{10})(a-\sqrt5)},
  \end{equation*}
  where
  $$ a:=\df{1+\sqrt5}{2}+\sqrt{\df{1+\sqrt5}{2}}.$$
  Lastly, an example from \cite[p.~772]{yi3}, with corrected sign errors, is given:
  \begin{equation}\label{99}
  \df{\varphi(e^{-\pi})}{\sqrt3\,\,\varphi(e^{-9\pi})}=2-\sqrt3-\df{\sqrt[3]{4}(5-3\sqrt3)}{(11\sqrt{3}-19)^{1/3}}
  -\left(2(11\sqrt3-19)\right)^{1/3}.
  \end{equation}
Comparing \eqref{99} with \eqref{999}, we see that the evaluations take rather different forms.

\section{An Incomplete Theta Function Evaluation;\linebreak Works of Seung Hwan Son and the Second Author}

On page 206 in his lost notebook \cite{lnb}, \cite[p.~180]{abII} Ramanujan recorded
the following identities.
\begin{entry} Let
\begin{equation}\label{p}
\dfrac{\varphi(q^{1/7})}{\varphi(q^7)} =1+u+v+w. 
\end{equation}
Then
\begin{equation}\label{p1}
p:=uvw=\dfrac{8q^2(-q;q^2)_{\infty}}{(-q^7;q^{14})_{\infty}^7}
 \end{equation}
  and
  \begin{equation}\label{p2}
\dfrac{\varphi^8(q)}{\varphi^8(q^7)}
-(2+5p)\dfrac{\varphi^4(q)}{\varphi^4(q^7)} +(1-p)^3=0.
\end{equation}
Furthermore,
\begin{equation}\label{p3}
u=\left(\dfrac{\alpha^2p}{\beta}\right)^{1/7},\quad
v=\left(\dfrac{\beta^2p}{\gamma}\right)^{1/7},\quad\text{and}
\quad w=\left(\dfrac{\gamma^2p}{\alpha}\right)^{1/7},
\end{equation}
 where $\alpha,\beta,$ and $\gamma$
are roots of the equation
\begin{equation}\label{p4}
r(\xi):=\xi^3+2\xi^2\left(1+3p-\dfrac{\varphi^4(q)}{\varphi^4(q^7)}
\right)+\xi p^2(p+4)-p^4=0. 
\end{equation}
For example,
\begin{equation}
\varphi(e^{-7\pi\sqrt7})=7^{-3/4}\varphi(e^{-\pi\sqrt7})
\left\{1+(-)^{2/7}+(-)^{2/7}+(-)^{2/7}\right\}. \label{ttt.1.6}
\end{equation}
\end{entry}

(We have corrected a misprint; Ramanujan wrote $7^{3/4}$ instead of $7^{-3/4}$ on the right-hand side of \eqref{ttt.1.6}.)  The terms  $(-)^{2/7}$  in \eqref{ttt.1.6} were not divulged by Ramanujan. To find the missing terms, one has to first solve the equation in \eqref{p2} for $\varphi^4(q)/\varphi^4(q^7)$ and choose the correct root. Second, one needs to solve the equation in \eqref{p4} for $\xi$ and use the roots in the correct order in \eqref{p3}. Ramanujan gave us no hints on how to do this. Since Ramanujan used the exponent $2/7$ on the right-hand side of \eqref{ttt.1.6}, we are almost certain that he had started to write down the same representation that we describe below. Our guess is that he stopped after finding $\alpha, \beta,$ and $\gamma$, but before he figured out their correct order.

The incomplete assertion \eqref{ttt.1.6} is one of only a few occasions in his notebooks where Ramanujan did not complete his formula.  We provide one example. On page 210 in his lost notebook, Ramanujan indicates that he has found the values of 14 specific  Rogers--Ramanujan continued fractions, but he gives the values of only three of them.  He evidently knew that he could perform all of  these evaluations, but he had other mathematical ideas that were more pressing to investigate. See \cite[pp.~62--75]{abI} for these evaluations.

In their book \cite{abII}, George Andrews and the first author were led by Ramanujan to correctly determine that in \eqref{p} \cite[p.~181]{abII}
 \begin{equation*}\label{u,v,w-def}
u:=2q^{1/7}\frac{f(q^5, q^9)}{\varphi(q^7)}, \qquad
v:=2q^{4/7}\frac{f(q^3, q^{11})}{\varphi(q^7)}, \qquad
w:=2q^{9/7}\frac{f(q, q^{13})}{\varphi(q^7)}.
\end{equation*}
 In a wonderful paper \cite{son}, Seung Hwan Son established proofs of \eqref{p1}--\eqref{p4}, which were reproduced in \cite[pp.~181--184]{abII}.  When Andrews and the first author wrote their second volume \cite{abII} on Ramanujan's lost notebook, they were unable to complete Ramanujan's evaluation of $\varphi(e^{-7\pi\sqrt7})/\varphi(e^{-\pi\sqrt7})$.

The completion of \eqref{ttt.1.6} has recently been accomplished by the second author \cite{rebak}; a brief sketch of his proof will now be given.

 As prescribed by Ramanujan in \eqref{ttt.1.6}, we set $q=e^{-\pi/\sqrt7}$ in \eqref{p}.  Next, from \eqref{p1}, $p=1$.  Then, we have to determine the correct root of the quadratic equation \eqref{p2}.  After doing so, we find that
\begin{equation*}
\frac{\varphi^4(q)}{\varphi^4(q^7)} = \frac{\varphi^4(e^{-\pi/\sqrt{7}})}{\varphi^4(e^{-\pi\sqrt{7}})} = 7.
\end{equation*}
We now have all the coefficients of the polynomial $r(\xi)$ defined in \eqref{p4}. We therefore need to determine the zeros $\alpha, \beta, \gamma$ of
\begin{equation*}
r(\xi) = \xi^3 - 6\xi^2 + 5\xi - 1,
\end{equation*}
which are
\begin{equation}\label{donedone}
\frac{1}{(2\cos\frac{k\pi}{7})^2}, \qquad k= 1, 2, 3.
\end{equation}

It would seem that  with \eqref{donedone}, we can determine $u, v, $ and $w$ in \eqref{p3}, and put their values in \eqref{p} to accomplished our goal of establishing the missing terms in \eqref{ttt.1.6}. However, it remains to determine the correct order of the roots $\alpha, \beta, \gamma$ in \eqref{p3}. The choice
\begin{equation*}
(\alpha, \beta, \gamma) = \Bigg(\frac{1}{(2\cos\frac{3\pi}{7})^2}, \frac{1}{(2\cos\frac{2\pi}{7})^2}, \frac{1}{(2\cos\frac{\pi}{7})^2}\Bigg)
\end{equation*}
is correct.  Hence,
\begin{equation*}
u = \Bigg(\frac{\cos\frac{2\pi}{7}}{2\cos^2\frac{3\pi}{7}}\Bigg)^{2/7}, \quad
v = \Bigg(\frac{\cos\frac{\pi}{7}}{2\cos^2\frac{2\pi}{7}}\Bigg)^{2/7}, \text{\quad and \quad}
w = \Bigg(\frac{\cos\frac{3\pi}{7}}{2\cos^2\frac{\pi}{7}}\Bigg)^{2/7}.
\end{equation*}
In conclusion, for $q=e^{-\pi/\sqrt7}$ the values of $u, v, $ and $w$ in \eqref{p} have been determined.  Ramanujan's incomplete formula \eqref{ttt.1.6} can now be made precise by
\begin{equation*}\label{done}
\varphi(e^{-7\pi\sqrt{7}}) = 7^{-3/4}\varphi(e^{-\pi\sqrt{7}})\Bigg\{1 + \bigg(\frac{\cos\frac{\pi}{7}}{2\cos^2\frac{2\pi}{7}}\bigg)^{2/7} + \bigg(\frac{\cos\frac{2\pi}{7}}{2\cos^2\frac{3\pi}{7}}\bigg)^{2/7} + \bigg(\frac{\cos\frac{3\pi}{7}}{2\cos^2\frac{\pi}{7}}\bigg)^{2/7}\Bigg\}\,.
\end{equation*}

In \cite{rebak}, the values of $\varphi(e^{-21\pi}), \varphi(e^{-35\pi}),$ and $\varphi(e^{-49\pi})$ are evaluated as well.

\section{Concluding Remarks}
Continuing our discussion above for determining further values of $\varphi(e^{-n\pi})$ from previously determined values, we can, in principal, apply the processes of  \emph{duplication} \eqref{dup}, \emph{dimidiation} \eqref{dim}, and \emph{change of sign} \eqref{sign}  to $\varphi(e^{-n\pi})$ to obtain values of $\varphi(\pm e^{-n2^a\pi})$, where $a\in \mathbb{Z}$. For example, we might attempt to use the values \eqref{3} and $\varphi(e^{-6\pi})$ given in \cite{yi}, to find values of $\varphi(e^{-3\cdot 2^a\pi})$.   However, because of the large number of manipulations, we would most likely need to invoke a computer denesting program.  Consequently, if the value of the requisite class invariant is known, it may be easier to directly apply the general procedures described above in an attempt to calculate a certain value of $\varphi(\pm e^{-n2^a\pi})$  in closed-form, which also may or may not be more elegant than what might be obtained by other means. Because of different approaches, different representations for the same $\varphi(e^{-n\pi})$ may arise,  as we demonstrated above for $\varphi(e^{-9\pi})/\varphi(e^{-\pi})$.

Readers will have observed that the methods of Yi in \cite{yi}, \cite{yi2}, and \cite{yi3} are apparently useful only when the degrees of the modular equations that are employed are `small.'  Likewise, the methods of the first author and Chan, and likely those of Ramanujan as well, become considerably more complicated to use for `larger' degrees.  In particular, we see from \eqref{Gn2} that class invariants for the square of the index $n$ in $\varphi(e^{-n\pi})$ are necessary. But the methods of all cited authors are similar.

The theta functions $\psi(q)$ and $f(-q)$ can be expressed in terms of $\varphi(q)$ \cite[Volume~2, p.~198]{nb}, \cite[pp.~39, 40, Entries 24, 25]{III}. Consequently, Ramanujan also established formulae for $\psi(q)$ and $f(-q)$ analogous to those for $\varphi(q)$, illustrated by \eqref{zalpha}  above \cite[Volume~2, pp.~210, 211]{nb}, \cite[pp.~123, 124, Entries 11, 12]{III}. See, for example, \cite[p.~326]{V} for several explicit values of $f(-q)$ recorded on page 250 in Ramanujan's first notebook \cite[Volume~1]{nb}. Ramanujan also evaluated a certain quotient of $\psi$-functions for several values of the parameters \cite[Volume~1, pp.~338, 339]{nb}, \cite[pp.~337--351]{V}. The ideas developed in \cite{yi} were extended by Yi and several others, and used to find similar values for $\psi(q)$, certain other products of theta functions, and the Rogers--Ramanujan continued fractions.

\subsubsection*{Acknowledgments}
The authors are grateful to the referee for a careful reading of their paper.

\end{document}